\documentclass[twoside, 12pt]{article}
\usepackage[latin1]{inputenc}
\usepackage{amscd}
\usepackage{amsmath}
\usepackage{amsfonts}
\usepackage{amssymb}
\usepackage{amsthm}
\usepackage{comment}
\usepackage[francais,english]{babel}
\usepackage[matrix,arrow]{xy}

\newcommand{\norm}[1]{\left\| #1 \right\|}

    \DeclareMathOperator{\Leb}{Leb}

\renewcommand{\geq}{\geqslant}
\renewcommand{\leq}{\leqslant}

\newcommand{\N}{\mathbb{N}}

\newcommand{\R}{\mathbb{R}}

\renewcommand{\phi}{\varphi}
\renewcommand{\epsilon}{\varepsilon}
\newcommand{\tq}{\ |\ }

\newcommand{\boR}{\mathcal{R}}
\newcommand{\boT}{\mathcal{T}}
\newcommand{\boQ}{\mathcal{Q}}

\newtheorem{thm}{Theorem}[section]
\newtheorem*{thm*}{Theorem}
\newtheorem{prop}[thm]{Proposition}
\newtheorem{lem}[thm]{Lemma}

\theoremstyle{definition}

\DeclareMathOperator{\lef}{left}
\DeclareMathOperator{\hor}{horz}
\newcommand{\km}{k_{\text{max}}}

\title{A note on stretched exponential decay of correlations for the
Viana-Alves map}

\author{V. Baladi\thanks{baladi@math.jussieu.fr, CNRS UMR 7586, 
Institut de Math\'e\-ma\-ti\-ques de Jussieu, 75005 Paris, France}
and S. Gou\"ezel\thanks{
sebastien.gouezel@ens.fr, Ecole Normale Supérieure,  75005 Paris, France}}

\date{November 2003}

\begin{document}

\maketitle

\section{Introduction}

Let $\phi:S^1\times \R$ to itself be given by
  \begin{equation}
  \phi(\omega,x)=(g(\omega),f(\omega,x))=
  (d \omega, a_0+\epsilon \sin(2\pi\omega)-x^2)
  \end{equation}
where $a_0\in (1,2)$ is fixed such that $x=0$ is a preperiodic point
for the map $h(x)=a_0-x^2$, and $d$ is an integer, say $\geq 16$.

For small $\epsilon>0$, this map leaves invariant a set of the form
$S^1\times I$ for some nonempty compact interval $I$.
It is known that this map has two positive
Lyapunov exponents Lebesgue-almost everywhere
(\cite{V}), that it has an ergodic SRB probability
(\cite{A}), and that the decay of the correlations
for this measure is faster than any polynomial
(\cite{ALP}). The aim of this work is to show that
the decay of the correlations is in fact at least
$O(e^{-c \sqrt{n}})$.

The main difference between our method and the method of
\cite{ALP} is that our construction is
inductive. In their article, if a point has many hyperbolic times
between $0$ and $N$ but has not yet been chosen, then it is not
in contradiction with Pliss' Lemma that this point does not have
hyperbolic times between $N$ and $2N$ for example. Thus, it is
possible that the measure of points remaining at time $2N$ is
quite large (and a careful study shows that, without new ideas,
their method will not give a decay rate better than $e^{-(\log
n)^2}$). In our inductive setting, everything restarts afresh
after each iteration, so we do not have this kind of problem.
This is made possible by a precise control of the geometry of the
system (while the result of \cite{ALP} is
valid in a much more general setting) -- in particular, we need
to use so-called \emph{hyperbolic returns} to control the size of
the sets given by the induction.

{\it This note was written in December 2002, when our result was
announced \cite{BG}. Since then, the second named author has
found the ``new ideas'' needed
to  enhance the  techniques in \cite{ALP},
obtaining a general abstract result \cite{G} which gives
as a particular case another proof of the present result.
We nevertheless believe that the ad hoc proof in this rough note
based on ideas from \cite{V, A, BBM}
should be made publicly available, at least on {\tt arxiv.org.}}

\section{Preliminary estimates}

We introduce a partition of $I$ (mod $0$) into the following
intervals:
  \begin{gather*}
  I_r=(\sqrt{\epsilon}e^{-r},\sqrt{\epsilon} e^{-(r-1)}) \text{ for
}r\geq 1,\\
  I_r=-I_{-r} \text{ for }r\leq -1,\\
  I_{0+}=I\cap [\sqrt{\epsilon},+\infty) \text{ and }I_{0-}=I\cap
(-\infty,-\sqrt{\epsilon}]
  \end{gather*}

We also write $I_r^+$ for the union of the three consecutive intervals
centered on $I_r$ (with the straightforward modifications for $I_{0+}$
and $I_{0-}$).

Given $(\omega,x)\in S^1\times I$, we define
$(\omega_j,x_j)=\phi^j(\omega,x)$.
Following \cite{V}, we take $\eta$ a positive
constant smaller than $1/3$ depending only on the quadratic map
$h$. We have (\cite[Lemma 2.1]{A})

\begin{lem}
\label{expanding_pres_critique}
There are constants $C_0,C_1>0$ such that for every small
$\epsilon>0$, we have an integer $N(\epsilon)$ satisfying
\begin{enumerate}
\item If $|x|<3\sqrt{\epsilon}$ then $\prod_{j=0}^{N(\epsilon)-1} |
\partial_x f(\omega_j,x_j)| \geq |x| \epsilon^{-1+\eta}$.
\item If $|x|<3\sqrt{\epsilon}$ then $|x_j|\geq \sqrt{\epsilon}$ for
$j=1,\ldots, N(\epsilon)$.
\end{enumerate}
\end{lem}

\begin{lem}
\label{expanding_hors_critique} There are $\sigma_2>1$ and $C_2>0$
such that $\prod_{j=0}^{k-1}\partial_x f(\omega_j,x_j) \geq C_2
\sigma_2^k$ whenever $|x_0|,\ldots,|x_{k-1}|\geq
e^{-9}\sqrt{\epsilon}$ and $|x_k|\leq 2\sqrt{\epsilon}$.

Moreover, $\prod_{j=0}^{k-1}\partial_x f(\omega_j,x_j) \geq C_2
\sqrt{\epsilon} \sigma_2^k$ whenever $|x_0|,\ldots,|x_{k-1}|\geq
e^{-9}\sqrt{\epsilon}$
\end{lem}

We say that the graph of a function $X:J\subset S^1 \to I$ (where $J$
is an interval) is an \emph{admissible curve} if $X$ is $C^2$ with
$|X'|\leq \epsilon$ and $|X''|\leq \epsilon$. Then
(\cite[Lemma 2.1]{V})
\begin{prop}
For small enough $\epsilon$, the image of an
admissible curve defined on an interval of length
$<1/d$ is still an admissible curve.
\end{prop}

\section{Construction of a Markov tower}

\subsection{Growing to a fixed size}

A \emph{rectangle} is a subset $R$ of $S^1\times I$ bounded by two
vertical lines, and two ``horizontal'' curves, i.e. graphs of
functions from a subset of $S^1$ to $I$. We shall write $\lef(R)$
for the left side of $R$, and $\hor(R)$ for the projection of $R$
on $S^1$. When the basis of the rectangle is $S^1$, it will always be
$(0,1)$, i.e. the possible discontinuity will always be at $0$.
A rectangle is \emph{admissible} if its horizontal
boundaries are admissible curves. An admissible rectangle is
\emph{gentle} if it is a subset of $S^1 \times
(\sqrt{\epsilon},1)$ containing $S^1 \times
(\sqrt{\epsilon},2 \sqrt{\epsilon})$, or the symmetric of such a
subset with respect to $S^1 \times \{0\}$, or if it contains
$S^1 \times (I_r\cup\ldots \cup I_{r+5})$ for some $r$ (in particular, 
the basis of a gentle rectangle will always be the full circle). 
This
definition is useful to avoid too thin rectangles and to keep
clean boundaries. In particular, the left boundary of a gentle
rectangle always contains an interval of the form $I_r$. 

We have the following analogue of Lemma 7.10 of \cite{BBM}:

\begin{lem}
\label{growing_to_fixed_size} There exist $q\in \N$ and $C>0$ such
that, for any admissible rectangle $R$ with basis $S^1$
and left boundary of size $\geq \epsilon^{1-\frac{3}{2}\eta}$,
there exists a partition $R_0,\ldots,R_s,R'_0,\ldots,R'_k$ of 
$R$ and times
$t_0,\ldots,t_s,t'_0,\ldots,t'_k \leq q$ such that:
\begin{enumerate}
\item For $0\leq i \leq s$, $\phi^{t_i}$ maps $R_i$ bijectively on 
$S^1 \times \Lambda$ for $\Lambda=I_1$ or $I_{-1}$, with distortion
bounded by $C$.
\item For $0\leq i\leq k$, the rectangle $\phi^{t'_i}(R'_i)$ is
admissible and gentle, its left boundary is of size $\geq
\epsilon^{1-\frac{3}{2}\eta}/C$, and the distortion of
$\phi^{t'_i}$ is bounded by $C$ on $R'_i$.
\item $\Leb(\bigcup R_i)\geq \Leb(R)/C$.
\end{enumerate}
\end{lem}
\begin{proof}
In this proof, every time we iterate the map, cut the rectangle
vertically in $d$, and apply the following procedure
independently to each part. Thus, at each step, the image of every
rectangle will have $S^1$ as its basis. From this point on, we will
only describe what happens in the $x$ direction.

Let $t$ be the first time
such that $\phi^t(R)$ meets $S^1\times \{|x|<e^{-9}
\sqrt{\epsilon}\}$.

If $\phi^t(R)$ also meets $S^1 \times \{|x|>3\sqrt{\epsilon}\}$,
then we can cut $\hor(\phi^t(R))\times \Lambda$ as a part of
$\phi^t(R)$, and hence subdivide $\phi^t(R)$ in three parts, for
which the return time will be $t$. This gives the required
construction: the number of iterates is bounded by a constant
$C(\epsilon)$ (according to the second part of Lemma
\ref{expanding_hors_critique}), the distortion is bounded since
in this finite number of iterates we have uniformly avoided the
critical point, and the vertical size is $\geq
(e^{-1}-e^{-2})\sqrt{\epsilon}$ at some point, whence it is $\geq
(e^{-1}-e^{-2})\sqrt{\epsilon} -2\epsilon$ on the left (because
the rectangle is bounded by admissible curves). This is $\geq
\epsilon^{1-\frac{3}{2}\eta}$ if $\epsilon$ is small enough.
Finally, the upper part $U$ will contain $\hor(U)\times
(\sqrt{\epsilon},2\sqrt{\epsilon})$, whence it is gentle, and the
lower part $V$ contains $\hor(V) \times (I_2\cup\ldots\cup I_7)$,
whence it is also gentle.

Otherwise, we set $(S_0,t_0)=(\phi^t(R),t)$. Note that $S_0
\subset S^1 \times \{|x|<3\sqrt{\epsilon}\}$. By Lemma
\ref{expanding_hors_critique}, $|\lef(S_0)|\geq C_2 |\lef(R)|\geq
C_2 \epsilon^{1-\frac{3}{2}\eta}$. We construct inductively
$(S_i,t_i)$ such that $S_i$ is a subset of $S^1 \times
\{|x|<3\sqrt{\epsilon}\}$, and with $|\lef(S_{i+1})|\geq C
\epsilon^{-\eta/2}|\lef(S_i)|$. This will imply that, if
$\epsilon$ is small enough, the process will stop after a finite
number $C(\epsilon)$ of iterates. Note that, with the process of
vertical cutting, $S_i$ will be replaced by a smaller $S'_i$, but with 
$|\lef(S'_i)|\geq |\lef(S_i)|-\epsilon \geq |\lef(S_i)|/2$ (since
$|\lef(S_i)|\geq C \epsilon^{1-\frac{3}{2}\eta})$), whence this will
only change the constants.

Assume $(S_i,t_i)$ is constructed. If $S_i$ meets $S^1\times
 \{ |x|<\epsilon^{1-\frac{3}{2}\eta}/(10e^9 C_2)\}$, we can cut a part of $S_i$ with a
horizontal line at height $\pm
\epsilon^{1-\frac{3}{2}\eta}/(5C_2)$ (recall that
$|\lef(S_i)|\geq \epsilon^{1-\frac{3}{2}\eta}/C_2$) and put it as
a $\phi^{t'_j}(R'_j)$ (for $t'_j=t_i$), such that the remaining part
$S$ satisfies $|\lef(S)|\geq |\lef(S_i)|/4$ and $S \subset
S^1\times \{|x|>\epsilon^{1-\frac{3}{2}\eta}/(5C_2)\}$. Note that
$\phi^{t'_j}(R'_j)$ will contain $\hor(\phi^{t'_j}(R'_j)) \times
(\epsilon^{1-\frac{3}{2}\eta}/(10e^8
C_2),\epsilon^{1-\frac{3}{2}\eta}/(10C_2))$ (there is a small loss
due to the fact that the boundaries are not straight lines). The
ratio of $e^8$ ensures that this interval contains at least $6$
consecutive $I_r$, and proves the gentleness of $\phi^{t'_j}(R'_j)$.
Moreover, it will satisfy $|\lef(\phi^{t'_j}(R'_j))|\geq
\epsilon^{1-\frac{3}{2}\eta}/(20C_2)$, which gives the claim on
its size.

Let $t$ be the first time such that $\phi^t(S)$ meets $S^1\times
\{|x|<e^{-9} \sqrt{\epsilon}\}$. If $\phi^t(S)$ also meets $S^1
\times \{|x|>3\sqrt{\epsilon}\}$, we cut it in three pieces as at
the beginning of the proof, and we stop the construction.
Otherwise, we set $(S_{i+1},t_{i+1})=(\phi^t(S),t_i+t)$. By Lemma
\ref{expanding_pres_critique}, we will have $t\geq N(\epsilon)$,
and during the first $N(\epsilon)$ iterates we will have an
expansion $\geq |x|\epsilon^{-1+\eta}\geq C \epsilon^{-\eta/2}$
(since $S \subset S^1\times
\{|x|>\epsilon^{1-\frac{3}{2}\eta}/5\}$). During the next
$t-N(\epsilon)$ iterates, we will have an expansion $\geq C_2
\sigma_2^{t-N(\epsilon)}$ according to Lemma
\ref{expanding_hors_critique}, which implies that $t\leq
C(\epsilon)$, and that globally the expansion will be at least
$C_2 C \epsilon^{-\eta/2}$. This proves the claim
$|\lef(S_{i+1})|\geq C \epsilon^{-\eta/2} |\lef(S_i)|$, and
concludes the construction.

We check that the desired properties are satisfied: the claims on
the size of the images come from the construction. The number of
steps in the construction is bounded, since at each step we have
an expansion of $C \epsilon^{-\eta/2}>1$. In each step, the number
of iterates is bounded by $C(\epsilon)$, thus the global number
of iterates is bounded. Finally, we iterate the map only outside
of the set $\{|x|<\epsilon^{1-\frac{3}{2}\eta}/(5C_2)\}$, which
implies that the distortion will be bounded. Finally, the claim on
$\Leb(\bigcup R_i)$ comes from the bounded distortion and the fact that
the number of rectangles will be bounded by $(2d)^q$.
\end{proof}

\subsection{Construction of the partition associated to an admissible
rectangle}

We fix $p_0=p_0(\epsilon)$ such that the expansion during a time
$p_0$ more than compensates for the distortion and the possible
contraction during the $q$ iterates of Lemma
\ref{growing_to_fixed_size}.

Write $r_j(\omega,x)=|r|$ if $x_j \in I_r$ with $|r|\geq 1$, $0$
otherwise. Consider $G_n(\omega,x)=\{1\leq i\leq n-1 \tq
r_i(\omega,x)\geq
\left(\frac{1}{2}-2\eta\right)\log\frac{1}{\epsilon}\}$. Take
$c>0$ small enough, and $c'>c$ very close to $c$. We say that $n$
is a hyperbolic return for $(\omega,x)$ if for every $0\leq k<n$,
we have
  \begin{equation*}
  \sum_{\substack{i\in G_n(\omega,x)\\k\leq i<n }} r_i(\omega,x) \leq
c'(n-k)
  \end{equation*}
and
  \begin{equation*}
  r_n(\omega,x)\geq 1.
  \end{equation*}

Write $H_n^*=\{(\omega,x) \tq n \text{ is the first hyperbolic
time }\geq p_0\}$. Then there exists $\gamma(\epsilon)>0$ and
$C(\epsilon)>0$ such that
  \begin{equation}
  \label{controle_retours_hyperboliques}
  \forall n\in \N,\ \Leb((S^1 \times I)-H_{p_0}^* \cup\ldots\cup H_{n-1}^*)\leq C e^{-\gamma
  \sqrt{n}}.
  \end{equation}

\begin{prop}
\label{proprietes_partition}
Let $R$ be a gentle admissible
rectangle. Then there exists a partition $\boR(R)=\bigcup_{n\geq
p_0}\boR_n(R)$ such that
\begin{enumerate}
\item $H_n^* \cap R \subset \bigcup_{k\leq n}\bigcup_{S\in \boR_k(R)} S$.
\item $\forall S\in \boR_n(R)$, the rectangle $f^n(S)$ has basis
$S^1$, satisfies
$|\lef(f^n(S))|\geq \epsilon^{1-\frac{3}{2}\eta}$,
and $f^n$ is uniformly expanding and has
uniformly bounded distortion on $S$.
\end{enumerate}
\end{prop}

\textbf{Construction of the initial partition $\boQ(R)$}

Since he wants a partition of the whole space, Alves starts from
the partition $\{ I_r \times S^1\}$. However, we start from an
admissible rectangle, whose boundary can have a slope $\epsilon$,
and in particular this boundary may cross $S^1 \times \{0\}$.
Thus, we have to construct a more complicated initial partition.

This partition $\boQ(R)=\{Q_i\}$ will have the following
properties:
\begin{enumerate}
\item Each $Q_i$ is an admissible rectangle,
contained in a set $S^1 \times I_r^+$, and its horizontal size is
of the form $1/d^{s}$ for some $s\in \N$.
\item $Q_i$ contains a set of the form $\hor(Q_i)\times I_r$.
\item If $Q_i \cap H_n^* \not=\emptyset$, then the horizontal size
of $Q_i$ is at least $1/d^{n}$.
\end{enumerate}
The last property is important because, if $Q_i$ intersects
$H_n^*$, we will try to iterate $Q_i$ exactly $n$ times, and we
need to recover a rectangle with basis $S^1$.

If the gentle rectangle $R$ is contained in $S^1 \times
\{x>\sqrt{\epsilon}\}$ and contains $S^1\times \{
\sqrt{\epsilon}<x <2\sqrt{\epsilon}\}$, it suffices to take
$Q_0=R$. So, we can assume that $R$ contains $\hor(R)\times
(I_a\cup\ldots \cup I_{a+5})$ for some $a$.

To construct $\boQ(R)$, we start from the partition $\boQ'$ in
sets $S^1 \times I_r$ for $r\leq (\frac{1}{2}-2\eta)\log(1/\epsilon)+2$
and $[\frac{k}{d^r},\frac{k+1}{d^r}] \times I_r$ for
$r>(\frac{1}{2}-2\eta)\log(1/\epsilon)+2$.

Note that if the horizontal boundary of $R$ intersects a $Q'$,
then it can intersect at most one rectangle of the same
horizontal size above or below (because of the bound $\epsilon$
on the slope and the smallness of the horizontal size).
Otherwise, for $r\leq (\frac{1}{2}-2\eta)\log(1/\epsilon)+2$, we would
have $\epsilon \geq |I_{r+1}|\geq C \epsilon^{1-2\eta}$, which is a
contradiction for $\epsilon$ small enough, and for 
$r> (\frac{1}{2}-2\eta)\log(1/\epsilon)+2$ we would have 
$|\hor(Q')| \epsilon \geq |I_{r+1}|$,
which implies $\epsilon d^{-r} \geq C \sqrt{\epsilon}e^{-r}$ and is again
a contradiction for $\epsilon$ small enough.
Thus, if we
form a block of three rectangles, one of them will be included in
$R$, and the intersection of this block with $R$ will give a
valid $Q_i$. This deals with the boundaries of $R$; in its
interior, simply put the remaining $Q'$ as $Q_j$. Note that the
gentleness of $R$ ensures that there will be no bad interaction
between the lower and upper boundaries of $R$.

We check the third claim on the hyperbolic returns: assume that
$Q_i \cap H_n^* \not=\emptyset$. If $|\hor(Q_i)|=S^1$, there is
nothing to prove. Otherwise, $Q_i=[\frac{k}{d^r},\frac{k+1}{d^r}]
\times I_r$ (or it comes from a block containing this), which
implies that its height is at most $e^{-r+2}$. Thus, for
$(\omega,x)\in Q_i$, $r(\omega,x) \geq r-2$ (because $r-2\geq
(1/2-2\eta)\log(1/\epsilon)$, we have $0\in
G_n(\omega,x)$), whence $\sum_{j\in
G_n(\omega,x)}r(\omega_j,x_j)\geq r-2$. Since $n$ is a hyperbolic
time for $(\omega,x)$, we obtain $r-2 \leq c'n$, whence $n\geq
r$ (as soon as $n>c'n+2$, which will be true for $\epsilon$ small
enough since $n>(r-2)/c'\geq \frac{1}{c'}(1/2-2\eta)\log(1/\epsilon)$).
Then $|\hor(Q_i)|=1/d^r \geq
1/d^{n}$.

\textbf{Construction of the partition $\boR(R)$}

Let $R$ be an admissible rectangle. An admissible subrectangle
$S$ of $R$ is $n$-\emph{good} if for every $0\leq j\leq n$, there
exists $r$ such that $\phi^j(\lef(S)) \subset I_r^+$, and there
exists $j\leq n$ such that $\phi^j(\lef(S))\supset I_r$, and $S
\cap H_n^*\not= \emptyset$.

Then there exists a partition
$\boR(R)=\bigcup_{n\geq p}\boR_n(R)$ such that
\begin{enumerate}
\item If $S \in \boR_n(R)$, then $\hor(S)$ is of the form
$\left[\frac{k}{d^{n}},\frac{k+1}{d^{n}} \right]$ for some
$0\leq k \leq d^{n}-1$.
\item $H_n^* \cap R \subset \bigcup_{k\leq n}\bigcup_{S\in \boR_k(R)} S$.
\item For every $0\leq j\leq n$ and $S \in \boR_n(R)$, there exists an
$I_{r_j}$ such that $\phi^j(\lef(S)) \subset I_{r_j}^+$.
\item For every $S\in \boR_n(R)$, either $S$ is $n$-good, or
there exists a $j$-good rectangle $T$
for some $j\leq n$, such that $S$ is subordinate to $T$.
(Our definition of subordinate is adapted from \cite{A}:
$S$ is subordinate to $T$ if, 
on the one hand there are $\ell \le j-1$
and $I_{r_\ell}$ with  $I_{r_\ell}\subset \varphi^\ell(\lef T)$,
and on the other hand $S$ is a subrectangle of
$\widetilde T$ with $\hor \widetilde T=\hor T$
and either the top admissible curve or the bottom admissible curve
of $\widetilde T$ coincides with that of $T$,
and
either $I_{r_\ell +1}$ or $I_{r_\ell-1}$ is included in
$\varphi^\ell(\lef \widetilde T)$.)
\end{enumerate}

In fact, it is sufficient to construct such a partition for each
$Q\in \boQ(R)$. And in this case, we can use more or less
directly the construction of Alves, up to checking that we have
enough control on horizontal sizes.

\textbf{Proof of proposition \ref{proprietes_partition}}

This is done in Alves, with minor modifications due to the fact
that our boundaries are not straight line but admissible curves.

The expansion with the hyperbolic returns, done in
\cite{BBM}, shows that the size at the end
will be $\geq \epsilon^{1-2\eta}/C$ for a constant $C$
independent of $\epsilon$. If $\epsilon$ is small enough, this
will be $\geq \epsilon^{1-\frac{3}{2}\eta}$. \qed

\subsection{Construction of the global partition}

We construct a partition $\boT$ of $X=S^1 \times \Lambda_\pm$ for
which the induced map will be Markov.

We start from the sets $K_\pm=S^1\times \Lambda_\pm$.
For each set  $S\in \boR_n(K_\pm)$, we apply Lemma
\ref{growing_to_fixed_size} to $\phi^n(S)$, giving 
$R_0,\ldots,R_s,R'_0,\ldots,R'_k$
and times $t_0,\ldots,t_s,t'_0,\ldots,t'_k$. 
Put the $n$-th preimage of each $R_i$ in the
partition $\boT$, with return time $n+t_i$. For $0\leq i\leq k$
apply inductively the construction process to
$\phi^{t'_i}(R'_i)$: decompose it as $\boR(\phi^{t'_i}R'_i)$, then use Lemma
\ref{growing_to_fixed_size} on each image, and go on.

Since, at each step, the process covers at least a proportion
$\geq 1/C$ of the remaining space (using Lemma
\ref{growing_to_fixed_size} and bounded distortion in between),
this will cover the whole space mod $0$. We write $R(\omega,x)$
for the return time of $(\omega,x)$ -- it is defined almost
everywhere.

\begin{thm}
There exists a constant $C$ such that $\Leb\{(\omega,x) \in X \tq
R(\omega,x)\geq  n\} \leq C e^{-\sqrt{n}/C}$.
\end{thm}

In the proof, we will use stopping time ideas, as introduced by
Young in \cite{Y}, but we will have to use
slightly different technical ideas, since the arguments of Young
would only give an estimate $C e^{-n^v}$ for every $v<1/2$ (this
is the decay rate she obtains when the estimate on the return
times is $e^{-\sqrt{n}}$). Note that the same technical idea can
be used to enhance her result, and we will indeed deduce from this
estimate on return times that the decorrelation rate is
$O(e^{-\sqrt{n}})$.

\begin{proof}
In the proof, we shall write $T_0(\omega,x) \leq
T_1(\omega,x)\leq \ldots T_{\km(\omega,x)}(\omega,x)$ for the
successive return times of $(\omega,x)$, i.e. the times given by
the use of Lemma \ref{growing_to_fixed_size}. While $(\omega,x)$
does not fall in a set $[k/d^q,(k+1)/d^q]\times \Lambda_\pm$ at
time $T_i$, then Proposition \ref{proprietes_partition} and Lemma
\ref{growing_to_fixed_size} give a next time $T_{i+1}$, and the
process stops only when the point returns to a set
$[k/d^q,(k+1)/d^q]\times \Lambda_\pm$, the return time being then
$T_{\km(\omega,x)}(\omega,x)=R(\omega,x)$.

Fix some $\delta>0$ very small. Since, at each step, a proportion
$1/C$ of the points returns, we have for each $n\in \N$
  \begin{equation}
  \label{estimee_exponentielle}
  \Leb\{ (\omega,x) \tq \km(\omega,x) \geq \delta \sqrt{n}\} \leq
  Ke^{-C(\delta) \sqrt{n}}.
  \end{equation}

Let $\tau_1<\ldots<\tau_i$ be fixed return time, and consider
  \begin{equation*}
  A(\tau_1,\ldots,\tau_i)
  =\{(\omega,x) \tq \km \geq i, T_1=\tau_1, T_2=\tau_2,\ldots, T_i=\tau_i\}.
  \end{equation*}
Let $R$ be a rectangle on which
$T_1=\tau_1,\ldots,T_{i-1}=\tau_{i-1}$, and write
$S=T^{\tau_{i-1}}(R)$. Then, by bounded distortion of
$f^{T_{i-1}}$ on $R$,
  \begin{equation*}
  \frac{ \Leb\{(\omega,x)\in R \tq T_i=\tau_i\} }{ \Leb(R)}
  \leq C \frac{\Leb\{(\omega,x)\in S \tq T_1=\tau_i-\tau_{i-1}\}}{\Leb S}.
  \end{equation*}
Lemma \ref{growing_to_fixed_size} gives that $\Leb(S) \geq
C(\epsilon)$, and
  \begin{equation*}
  \{(\omega,x)\in S \tq T_1=\tau_i-\tau_{i-1}\}
  \subset S^1 \times I - (H_{p_0}^*\cup \ldots \cup
  H_{\tau_i-\tau_{i-1}-1}^*)
  \end{equation*}
whence by Equation \eqref{controle_retours_hyperboliques},
  \begin{equation*}
  \Leb\{(\omega,x)\in S \tq T_1=\tau_i-\tau_{i-1}\} \leq C
  e^{-\gamma \sqrt{\tau_i-\tau_{i-1}}}.
  \end{equation*}
Summing these equations on all rectangles $R$, we obtain
  \begin{equation*}
  \Leb(A(\tau_1,\ldots,\tau_i)) \leq C
  \Leb(A(\tau_1,\ldots,\tau_{i-1}))e^{-\gamma
  \sqrt{\tau_i-\tau_{i-1}}}.
  \end{equation*}
Let us write $a_n=C e^{-\gamma \sqrt{n}}$. We get
  \begin{equation*}
  \Leb(A(\tau_1,\ldots,\tau_i)) \leq a_{\tau_1}
  a_{\tau_2-\tau_1}\ldots a_{\tau_i-\tau_{i-1}}.
  \end{equation*}
Summing finally on all possible sequences $\tau_1<\ldots<\tau_i$
with $\tau_i \geq n$ and $i\leq \delta \sqrt{n}$, we get
  \begin{equation}
  \label{on_avance_sur_Leb}
  \begin{split}
  \Leb\{(\omega,x)\tq \km(\omega,x)<\delta \sqrt{n}, R(\omega,x) \geq n\}
  &\leq \sum_{i\leq \delta \sqrt{n}} \sum_{\substack{\tau_1<\ldots<\tau_i\\
  \tau_i\geq n}} a_{\tau_1}
  a_{\tau_2-\tau_1}\ldots a_{\tau_i-\tau_{i-1}}
  \\&
  =\sum_{i\leq \delta \sqrt{n}} \sum_{\substack{j_1+\ldots+j_i
  \geq n\\j_1,\ldots,j_i>0}}
  a_{j_1}\ldots a_{j_i}.
  \end{split}
  \end{equation}

\begin{lem}
\label{lemme_technique_stretched}
Let $a_n$ be a sequence with
$a_n=O(e^{-\gamma \sqrt{n}})$. Then there exists $D>0$ such that
$b_n=D^{-1} a_n$ satisfies
  \begin{equation*}
  u_n:=\sum_{i=0}^n \sum_{\substack{j_1+\ldots+j_i=n\\j_1,\ldots,j_i>0}} b_{j_1}\ldots b_{j_i}
  =O(n^2 e^{-\gamma \sqrt{n}}).
  \end{equation*}
\end{lem}
Let $D$ be given by the lemma. Then Equation
\eqref{on_avance_sur_Leb} implies that
  \begin{equation}
  \label{on_arrive_sur_Leb}
  \begin{split}
  \Leb\{(\omega,x)\tq \km(\omega,x)<\delta \sqrt{n}, R(\omega,x) \geq
  n\}&
  \leq
  D^{\delta \sqrt{n}}\sum_{i\leq \delta \sqrt{n}} \sum_{p=n}^\infty \sum_{\substack{
  j_1+\ldots+j_i=p\\ j_1,\ldots,j_i>0}}
  b_{j_1}\ldots b_{j_i}
  \\&
  \leq D^{\delta\sqrt{n}} \sum_{p=n}^\infty u_p
  \leq D^{\delta \sqrt{n}} E n^3 e^{-\gamma \sqrt{n}}.
  \end{split}
  \end{equation}
We now choose $\delta$ small enough so that $D^\delta
e^{-\gamma}<1$ (note that $D$ does not depend on $\delta$), and
Equation \eqref{on_arrive_sur_Leb} gives a bound of the form $S
e^{-\gamma' \sqrt{n}}$. Add finally Equation
\eqref{estimee_exponentielle}, to get
  \begin{equation*}
  \Leb\{(\omega,x)\tq R(\omega,x) \geq
  n\}\leq T e^{-\min(\delta,\gamma') \sqrt{n}}
  \end{equation*}
which is the conclusion of the theorem.
\end{proof}

\begin{proof}[Proof of Lemma \ref{lemme_technique_stretched}]
Writing $s\star t$ for the convolution of the sequences $s_n$ and
$t_n$, i.e. $(s\star t)_n=\sum_{i=0}^n s_k t_{n-k}$, then
$u=\sum_{j=0}^\infty b^{\star j}$.

Write $w_n$ for a sequence equal to $e^{\gamma \sqrt{n}}/n^2$ for
$n$ large enough, and satisfying $w_{n+p}\leq w_n w_p$. Let
$\norm{s}=\sum w_n s_n$ for a sequence $s_n$ such that this sum
is finite. Then $\norm{s\star t}\leq \norm{s}\norm{t}$.

In particular, $\norm{a}<\infty$, whence for $D$ large enough,
$\norm{b}=\norm{a}/D<1$. Then $\norm{b^{\star j}}\leq \norm{b}^j$,
and $\norm{u} \leq \sum \norm{b}^j=\frac{1}{1-\norm{b}}<\infty$.
In particular, $w_n u_n$ is bounded, i.e. $u_n=O(n^2 e^{-\gamma
\sqrt{n}})$.
\end{proof}

\section{Decay of correlations}

It is not difficult to obtain an aperiodic tower.
Then, the rate of decay of correlations may be obtained by using
the coupling argument in \cite{Y} combined with our
technical lemma 
\ref{lemme_technique_stretched}.

\bibliographystyle{alpha}

\end{document}